\documentclass[11pt,a4paper]{article}

\title{Navier-Stokes and stochastic Navier-Stokes equations via Lagrange multipliers}
\author{Ana Bela Cruzeiro
\vspace{3mm}\\
{\footnotesize GFMUL and Dep. de Matem\'atica  
Instituto Superior T\'ecnico,}\\
{\footnotesize  Av. Rovisco Pais,
1049-001 Lisboa, Portugal }}
\date{}

\usepackage{amssymb,amsmath,amsfonts,amsthm,color}

\usepackage{framed}

\theoremstyle{plain}

\theoremstyle{definition}

\catcode`\@=11
\def\mequal{\mathrel{\mathpalette\@mvereq{\hbox{\sevenrm m}}}} 
\def\@mvereq#1#2{\lower.5\p@\vbox{\baselineskip\z@skip\lineskip1.5\p@
    \ialign{$\m@th#1\hfil##\hfil$\crcr#2\crcr=\crcr}}}
\def\partr#1#2{/\kern-.08333em/_{#1,#2}^{\phantom{.}}}
\def\invpartr#1#2{/\kern-.08333em/_{#1,#2}^{-1}} 
\def\hpartr#1#2{/\kern-.08333em/_{#1,#2}^{h}}
\def\Epartr#1#2{/\kern-.08333em/_{#1,#2}^{E}}

\def\newdot{{\kern.8pt\cdot\kern.8pt}}
\def\,{\relax\ifmmode\mskip\thinmuskip\else\thinspace\fi}
\def\{{\relax\ifmmode\lbrace\else $\lbrace$\fi}
\def\}{\relax\ifmmode\rbrace\else $\rbrace$\fi}

\def\L{\mathcal L}
\def\P{\mathbb{P}}
\def\H{\mathbb{H}}

\def\T{\mathbb{T}}

\def\div{\textup{div}}

\font\sevenrm=cmr7

\renewcommand{\theequation}{\arabic{section}.\arabic{equation}}

\begin{document}
\maketitle
\makeatletter 
\renewcommand\theequation{\thesection.\arabic{equation}}
\@addtoreset{equation}{section}
\makeatother 

\begin{abstract}

\it  We show that  the Navier-Stokes  as well as a random perturbation of this equation can be derived from a stochastic variational principle 
where the pressure is introduced as a Lagrange multiplier. Moreover we describe how to obtain corresponding constants
of the motion.\rm

\end{abstract}

\vskip 7mm

\section{ \bf Introduction}\label{section1.}

\quad  Navier-Stokes equation describes the velocity of incompressible viscous fluids. We  consider this equation
with periodic boundary conditions; namely, if $v=v(t,x), ~ x\in \T$, denotes this velocity at time $t$, with $\T$ being the
$d$-dimensional flat torus that we identify with $[0,2\pi ]^d$, it reads

\begin{equation}\label{eq1.1}
\frac{\partial v}{\partial t}+ (v\cdot\nabla)v =\nu \Delta v-\nabla p,\quad \div(v )=0,
\end{equation}
where $\nu$ is a positive constant (the viscosity coefficient) and $t\in [0,T]$. The function $p=p(t,x)$ denotes
the pressure and is also an unknown in the equation.

\quad Lagrange's point of view consists in  describing  positions of particles:
 it concerns the flows driven by the velocity fields. Lagrangian trajectories for the Euler equation (the case where there
 is no viscosity term) have been identified as minimisers of the kinetic energy defined on the space of diffeomorphisms by
 V. Arnold in \cite{Arnold}. In other words they are geodesics for a $L^2$ metric on such space of curves. This geometric approach
 to the Euler equation was developed in the fundamental paper by D. Ebin and J. Marsden (\cite{EM}) and gave rise to many subsequent
 works. It is well known that the pressure in an incompressible fluid acts like a Lagrange multiplier and one can, indeed, derive
 the Euler equation from a variational principle with such a multiplier (cf. for example \cite{Cons, H, S}).
 
\quad  Navier-Stokes equation, being a dissipative physical system, does not correspond to  analogous deterministic variational principles. 
 Nevertheless, by replacing the Lagrangian flows by stochastic ones, we may still derive this equation from a (stochastic) variational
 principle associated with the energy. Then  the velocity field is identified with the drift of the Lagrangian diffusion process,  which is a time derivative
 after conditional expectation of the paths. Inspired by
  \cite{NYZ} and \cite{Y}, such a stochastic variational principle was proved in \cite{CC}. More recently it was generalised in the context of Lie groups in 
  \cite{ACC} and many other dissipative systems can be derived using the same kind of ideas (cf. also \cite{CCR}). Moreover stochastic partial differential
   equations were also obtained by variational principles, corresponding to random perturbations of the action functionals, in \cite{CCR}. We refer 
   to \cite{H1} and other subsequent works from the same author, where a  different variational approach to stochastic fluid dynamics is
   developed (to derive stochastic partial differential equations).
  
\quad  In this paper we show that it is  possible to derive the Navier-Stokes equation from a (stochastic) variational principle with a Lagrange
  multiplier  expressed in terms of  the pressure. Although we consider here a flat case, the principle can be extended to general manifolds following the
  construction in \cite{AC}. For the general theory of stochastic differential equations on manifolds we refer for example to \cite{IW}.
  
 \quad Stochastic Noether's theorem was introduced in \cite{TZ1}, \cite{TZ2}  in the context of stochastic processes associated with the heat equation. A conserved
  quantity corresponds there to a martingale. In the spirit of this theorem as well as of \cite{CL}, we present  a  result about conserved quantities associated to our stochastic variational principle. The main difference with the the one  of 
\cite{CL} is that we consider here the Lagrangian motion as a stochastic  flow (with respect to its initial values $x$) and in the notion of symmetry we integrate
with respect to the variable  $x$. 

\quad It should be stressed, here, that in our derivation of the Navier-Stokes equation no random perturbation is added. What we advocate is an approach where the presence
of the Laplacian in Navier-Stokes equation is interpreted as the underlying presence of diffusion processes, used afterwards for studying (1.1) in probabilistic terms.
In the last section we show how to derive a variational approach to a randomly perturbed Navier-Stokes equation ((4.3)).

 \vskip 7mm
  
\section{ \bf A stochastic variational principle}\label{section2.}
 
\quad  On a fixed standard probability space $(\Omega, \P, P)$ endowed  with an increasing filtration $\P_t$ that satisfies the standard assumptions,
we consider $\xi$ to be a semimartingale with values in $\T$,  namely
\begin{equation}
d\xi_t (x)=dM_t (x) +D_t \xi (x) dt,\qquad \xi_0 (x)=x,
\end{equation}
where $x\in \T$, $M_t$ is the martingale part in the decomposition of $\xi_t$  and $D_t \xi$ its drift (for simplicity we do not write the probability parameter $\omega \in \Omega$ 
in the formulae).

\quad Recall the definition of generalised derivative, that we denote by $D_t$: for
$F$ defined in $[0,T] \times \T $,

\begin{equation}
D_t F(t, \xi_t (x))=\lim_{\epsilon \rightarrow 0} \frac{1}{\epsilon}E_t [F(t+\epsilon , \xi_{t+\epsilon}(x)) - F(t ,\xi_t (x))]
\end{equation}
when such (a.s.) limit exists, and where $E_t$ denotes the conditional expectation with respect to $\P_t$. This definition justifies in particular
the notation used in (2.1), since the generalised derivative corresponds to the seminartingale's drift.

\quad If $W_t$  is a $\P_t$-adapted Wiener process, we denote by
 $g_t (\cdot )$ diffusions on the torus $\T$ of the form

\begin{equation}
dg_t (x)=\sqrt{2\nu} dW_t +v(t, g_t (x))dt,\qquad g_0 (x)=x 
\end{equation}
with $x\in \T$, $dW_t $  the It\^o differential. The drift function $v$ is assumed to be regular enough so that $g_t (\cdot )$ are
diffeomorphisms (cf. \cite{K}).

\quad Note that we do not require, a priori, the vector field $v$ to be divergence free.

\quad For the particular cases $F(t,x) =x$ and $F(t,x) =v(t,x)$, we have, respectively,

$$ D_t g_t (x)= v(t, g_t (x))$$
and, using It\^o's formula,
\begin{equation}
 D_t v(t,g_t (x))=D_t D_t g_t (x)=\big( \frac{\partial}{\partial t} v +(v.\nabla )v +\nu \Delta v \big) (t, g_t (x)).
\end{equation}

\vskip 5mm

Let $\H$ be a linear subspace dense in $L^2 ([0,T]\times \T)$. Define the action functional

\begin{align}
S(g,p)&=\frac{1}{2} E\int_0^T \int |D_t g_t (x)|^2 dt dx +E\int_0^T \int p(t, g_t (x)) (\det \nabla g_t (x) -1)dtdx \\
&:=S^1 (g,p) +S^2 (g,p)
\end{align}
for $p\in \H$   and where $E$ denotes expectation (with respect to $P$).

\vskip 5mm
We consider  variations

$$
g_t (\cdot )\rightarrow g_t^\epsilon (\cdot ) =g_t (\cdot ) +\epsilon  h (t,g_t (\cdot ))
$$

$$
p(t, \cdot )\rightarrow p^\epsilon (t,\cdot )= p(t, \cdot ) +\epsilon \varphi (t,g_t (\cdot ))
$$
with $h(t,x)$ and $\varphi (t,x )$ deterministic and smooth in $x$, $\varphi \in \H$.
We also assume that $h(T,\cdot )=h(0,\cdot )=0$.
\vskip 5mm
\quad We have  the following
\vskip 3mm

\bf Theorem. \rm A diffusion $g_t $   of the form (2.3) and a function $p\in \H$ are critical for the action functional (2.5) iff the drift 
$v(t,\cdot )$ of $g_t$ satisfies
the Navier-Stokes equation (without external force)

\begin{equation}
\label{NS}
\partial_ t v +(v.\nabla )v =\nu \Delta v -\nabla p,\qquad \hbox{div} ~v(t,\cdot)=0,
\end{equation}
with $x\in \T, t\in [0,T]$.
\vskip 5mm
\bf Proof. \rm Using the notation $\delta S (g,p)= \left.\frac{d}{d\varepsilon}\right|_{\varepsilon=0} S(g^\epsilon ,p^\epsilon )$, the 
variation of the   first term in the action gives:

$$\delta S^1 (g,p)=E\int_0^T \int (D_t g_t (x).D_t h (t,g_t (x)))dt dx.$$

The notation $<\cdot ,\cdot >$ stands below for the $L^2 (\T )$ scalar product.
By It\^o's formula, the expression
$$
d<D_t g_t ,h> - <D_t D_t g_t ,h> dt-<D_t g_t , D_t h>dt- <dDg_t , dh>
$$
where the last term denotes the It\^o contraction, is the differential of a  martingale (whose expectation vanishes); therefore
\begin{equation} 
D_t < D_t g_t . h> =<D_t D_t g_t .h>+<D_t g_t ,D_t h> +<dD_t g_t , dh>.
\end{equation}
We deduce that
$$\delta S^1 =E <Dg_T,h(T, g_T )>-E<Dg_0,h(0, g_0 )> - E\int_0^T \int (D_t D_t g_t (x). h(t,g_t (x) ))dtdx$$

$$-E\int_0^T \int (dD_t g_t (x). dh(t, g_t (x)))  dx.$$

$$=- E\int_0^T \int (D_t D_t g_t (x). h(t,g_t (x) ))dtdx -2\nu E\int_0^T (\nabla v. \nabla h) (t, g_t (x))dt dx$$

$$= -E \int_0^T \int ((\partial_t v +(v.\nabla )v-\nu \Delta v). h)(t, g_t (x)$$
where, for the last equality we have used the equality $D_t D_t g_t (x) =D_t v(t, g_t (x))= (\partial_t v +(v.\nabla )v+\nu \Delta v)(t, g_t (x)$
and integration by parts.

Concerning the second part of the action functional, we have

\begin{align}
\delta S^2 &= E\int_0^T \int \varphi (t, g_t (x))  (\det \nabla g_t (x) -1)dtdx \\
& + E\int_0^T \int (\nabla p (t,g_t (x)) .h(t,g_t (x)) (\det \nabla g_t (x) -1)dtdx\\
&+ E\int_0^T \int p(t,g_t (x)) \left.\frac{d}{d\varepsilon}\right|_{\varepsilon=0} \det \nabla (g_t (x)+\epsilon h(t,g_t (x)) dt dx
\end{align}

\vskip 5mm

Since $\varphi$ is arbitrary we conclude from (2.9) that critical points of the action are volume-preserving diffeomorphisms
($\det \nabla g_t (x) =1$)
and therefore have divergence-free drifts. It follows immediately that $(2.10)=0$ so we only have to compute (2.11). We have,

$$ \left.\frac{d}{d\varepsilon}\right|_{\varepsilon=0} \det \nabla (g_t (x)+\epsilon h(t,g_t (x))) =\det \nabla g_t (x)~ \hbox{tr}
\Big( (\nabla g_t (x))^{-1}  \left.\frac{d}{d\varepsilon}\right|_{\varepsilon=0}\nabla g_t^\epsilon (x) \Big)$$

$$=\det (\nabla g_t (x)) ~\hbox{tr} \big( (\nabla g_t (x))^{-1} \nabla (h(t,g_t (x))) \big).$$

Since

$$\partial_i (p(t,g_t  ) (\nabla g_t )^{-1}_{ij} h(t,g_t )^j )= \partial_i (p(t,g_t ))(\nabla g_t )^{-1}_{ij} h(t,g_t )^j +p(t,g_t )(\nabla g_t )^{-1}_{ij} \partial_i (h^j (t,g_t  ))$$
$$+p(t,g_t ) h^j (t,g_t ) \partial_i (\nabla g_t )^{-1}_{ij}, $$
and we are in the periodic case,
$$(2.11)= -E\int_0^T \int [\partial_i (p(t,g_t ))(\nabla g_t )^{-1}_{ij} +p(t,g_t )\partial_i ((\nabla g_t )^{-1}_{ij} )] h^j (t,g_t  ) \det \nabla g_t ~dtdx.$$
Notice that we already concluded that $\det \nabla g_t =1$. On the other hand,
$$\sum_i \partial_i  (\nabla g_t )^{-1}_{ij} =0.$$
Indeed, derivating the equality
$\det \nabla g_t =1$, we get

$$\partial_k \det (\nabla g_t )= ~\hbox{tr}\big( (\nabla g_t )^{-1}\partial_k (\nabla g_t )\big)
=\sum_i (\nabla g_t )^{-1}_{ij}\partial_k \partial_i g_t^j =0.$$
Also, derivating equality

$$(\nabla g_t )^{-1}_{ij}\partial_k g_t^j =\delta_{ik}$$
we obtain
$$\sum_i \partial_i (\nabla g_t )^{-1}_{ij}\partial_k g_t^j +(\nabla g_t )^{-1}_{ij}\partial_i \partial_k g_t^j =0;$$
therefore
$$\sum_i \partial_i (\nabla g_t )^{-1}_{ik}=-\big( (\nabla g_t )^{-1}_{ij}\partial_k \partial_i g_t^j \big)(\nabla g_t )^{-1}_{jk}=0$$
and 
$$(2.11)=  -E\int_0^T \int (\partial_i (p(t,g_t  (x)))(\nabla g_t (x))^{-1}_{ij}) h_t^j (g_t (x)) \det \nabla g_t (x) )dtdx $$
$$=-E\int_0^T \int (\nabla p (t,g_t  (x)).h(t,g_t (x)))dtdx.$$

Putting together the expressions for
$\delta S^1$ and $\delta S^2$, we conclude that $\delta S=0$ in the class of variations considered,  is equivalent to the condition 
$$E \int_0^T \int (\partial_ t v +(v.\nabla )v -\nu \Delta v+\nabla p (t,g_t (x)))h(t, g_t (x)) dt dx=0$$
 for every test function
$h$, together with the incompressibility condition $\det \nabla g_t (x)=1$.

\vskip 7mm

\bf Remark 1. \rm Comparing with \cite{CC, ACC}, the variations we have used here are defined by shifts, since we do not have to work a priori in the class of measure-preserving flows.

\vskip 2mm

\bf Remark 2. \rm It is possible to consider subspaces $\H$ which are not dense in $L^2 ([0,T]\times \T)$. In this case the resulting equation of motion is the projection of the Navier-Stokes one in the corresponding space.

\vskip 7mm
\section{ \bf On  conserved quantities}\label{section3.}

\quad In this section we present a Noether-type result where only transformations in space of the Lagrangian function are considered. 
A more  general study of symmetries for  equations obtained by stochastic variational principles will be considered in a forthcoming work.

\quad Let us  consider transformations of the following form:
$$
g_t (\cdot )\rightarrow g_t^\alpha (\cdot ) =g_t (\cdot ) +\alpha  \eta (t, g_t (\cdot ))
$$
with $\eta$ smooth, $\eta (0,\cdot )=\eta (T,\cdot )=0$. We say that the Lagrangian

\begin{align}
L (g,p) &= \frac{1}{2}|D_t g_t (x)|^2  + p(t, g_t (x)) (\det \nabla g_t (x) -1)\\
&:=L^1 (g,p) +L^2 (g,p)
\end{align}
used in the definition of the action functional (2.5), is invariant under the transformation associated with $\eta$ if there
exits a  function  $G: [0,T]\times \T \rightarrow \mathbb R$ such that for every $t$, $P$-a.e.,

$$\left.\frac{d}{d\alpha}\right|_{\alpha =0}
\int L(g_t^{\alpha}  ,p) dx= \int D_t G (t, g_t (x) )dx.$$

\vskip 3mm

\bf Theorem. \rm If $L$ is invariant under the transformation associated with $\eta$ then,
denoting $\L_t = \frac{\partial}{\partial t} +(v\cdot \nabla ) +\nu \Delta $ 
where $v(t,\cdot )$ is the solution of the Navier-Stokes equation considered above,
the following identity 
$$\int \Big( \L_t (v \eta -G )\Big) (t,x) dx =0 $$
holds for all $t\in [0,T]$.

\vskip 5mm
\bf Proof. \rm 
Considering the first term in the Lagrangian, we have
$$\left.\frac{d}{d\alpha}\right|_{\alpha =0} L^1 (g_t^{\alpha} ,p) =(D_t g_t (x). D_t \eta (t, g_t (x)))$$
and, by the arguments in the proof of last section's Theorem,

$$
\left.\frac{d}{d\alpha}\right|_{\alpha =0} L^2 (g_t^{\alpha} ,p) =
(\nabla p (t, g_t (x))\cdot \eta (t, g_t (x)) (det \nabla g_t (x)-1) 
$$
\begin{equation}
+\partial_i (p (t, g_t (x))\nabla g_t (x)^{-1} \eta^j (t, g_t (x)))-(\nabla p (t, g_t (x)). \eta (t, g_t (x))).
\end{equation}
We know that $ (\det \nabla g_t (x)-1)=0$ on the critical points of the action functional, therefore the first
term in the r.h.s. of last equality vanishes. The second one also vanishes after integration in $x$, as we
consider periodic boundary conditions. We are therefore left with the equality, valid for $g_t$ critical of the
action functional,

$$\int \Big( (D_t g_t (x).D_t \eta (t, g_t (x))) 
- (\nabla p (t, g_t (x)). \eta (t, g_t (x)))\Big) dx = \int D_t G (t, g_t (x) )dx,$$
$P$-a.e. Using the identity

$$
D_t ( D_t g_t (x).\eta (t, g_t (x))) =(D_t D_t g_t (x).\eta (t, g_t (x))) +(D_t g_t (x). D_t \eta (t, g_t (x)))$$
\begin{equation}
+(dD_t g_t (x). d\eta (t, g_t (x)).
\end{equation}
 From the two last equalities 
we deduce that

$$\int D_t \Big( ( D_t g_t (x). \eta (t, g_t (x))
 - G (t, g_t (x) )\Big) dx =0.$$
 We have $D_t g_t (x)=v(t, g_t (x))$.
By the incompressibility condition, the flow $g_t (\cdot )$ keeps the measure $dx$ invariant (a.s.) and the result follows
from the expression of the operator $D_t$.

\vskip 3mm

\quad Comparing with the finite-dimensional Noether's theorem of  \cite{TZ1}, \cite{TZ2}, here we have an extra integration with respect to the space variable $x$
in the derived notion of conserved quantities.

\vskip 7mm
\section{ \bf A stochastic Navier-Stokes equation}\label{section4.}
In this section we show that it is also possible to derive random perturbations of the Navier-Stokes equation from a stochastic variational principle.

\quad Let $\xi$ be a semimartingale with values in $\T$ of the form (2.1).
We consider  the random action functional

$$
\tilde S(\xi ,p)=\frac{1}{2} \int_0^T \int |D_t \xi_t (x)|^2 dt dx +\int_0^T \int (D_t \xi_t (x). dM_t (x))
- \sqrt{2\nu} \int_0^T \int (D_t \xi_t (x). dW_t  )
$$
\begin{equation}
+\int p(t, \xi_t (x)) (\det \nabla \xi_t (x) -1)dtdx, 
\end{equation}
with $p \in \H \subset L^2 ([0,T]\times \T)$.
Variations of $\xi $ and $p$ are taken as in Section 1, namely

$$
g_t (\cdot )\rightarrow g_t^\epsilon (\cdot ) =g_t (\cdot ) +\epsilon  h (t,g_t (\cdot ))
$$

$$
p(t, \cdot )\rightarrow p^\epsilon (t,\cdot )= p(t, \cdot ) +\epsilon \varphi (t,g_t (\cdot ))
$$
except that here we allow $h$ and $\varphi$ to be random.

\quad We want to characterise critical points of $\tilde S$ of the form
$$dg_t (x)=\sqrt{2\nu} dW_t +v(t, g_t (x))dt,\qquad g_0 (x)=x$$
now considering  the vector field $v$ to be random. We proceed as in the theorem of section 1. The computations are analogous and
we have to add, in the
 variations of $S$, those of
the second and third new terms of this functional. These terms give,

$$\int_0^T \int [(h(t, g_t ). \sqrt{2\nu} dW_t )+ (D_t g_t .(\nabla h (t, g_t ).dW_t ))-(h(t, g_t ). \sqrt{2\nu} dW_t )]dx$$
that reduces to
$$\int_0^T \int v(t, g_t (x). (\nabla h (t, g_t (x)).dW_t )dx= \int_0^T v(t,x).(\nabla h(t,x).dW_t )$$
$$=- \int_0^T ((\nabla v(t,x).h(t,x)).dW_t )$$
equality which holds $P$-almost surely.

\quad We therefore conclude that a diffusion process of the form 
$$dg_t (x)=\sqrt{2\nu} dW_t +v(t, g_t (x))dt, \quad g_0 (x)=x$$
is critical for the action functional $\tilde S$ iff its (random) drift $v(t,\cdot )$ satisfies the following Navier-Stokes
stochastic partial differential equation:
\begin{equation}
d v +(v.\nabla )v =\sqrt{2\nu} \nabla v .dW_t + \nu \Delta v -\nabla p,\qquad \div ~v(t,\cdot)=0,
\end{equation}

with $x\in \T, t\in [0,T]$.

\quad This stochastic equation can be also regarded as a (Stratonovich) perturbation of the Euler one.
 Indeed, denoting by $\circ dW$ the Stratonovich differential, it can be written as

\begin{equation}
d v +(v.\nabla )v =\sqrt{2\nu} \nabla v \circ dW_t -\nabla p, \qquad \div ~v(t,\cdot)=0.
\end{equation}

\vskip 10mm

\noindent{\bf Acknowledgements}\\
The author was  supported by FCT Portuguese grant PTDC/MAT-STA/0975/2014.

She wishes to thank the anonymous referees for a careful reading of the first manuscript of this paper.

\end{document}